\documentclass[letterpaper, 10 pt, conference]{ieeeconf}  

\usepackage{url}  
\usepackage{graphicx}  
\usepackage{amsmath,amssymb,amsfonts,mathrsfs,bm} 
\usepackage{caption}
\usepackage{cleveref}
\usepackage{array}
\usepackage{scalerel}

\usepackage{algorithm}
\usepackage{algpseudocode}
\usepackage{booktabs}
\usepackage{epstopdf}
\usepackage{multirow}
\usepackage{rotating}
\usepackage{algpseudocode}
\usepackage{cite}
\usepackage{epsfig}
\usepackage{tabularx}
\usepackage[table]{xcolor}
\usepackage{epstopdf}
\usepackage{accents}
\usepackage{mathtools}
\usepackage{epstopdf}
\usepackage{pgfplots}

\newcolumntype{P}[1]{>{\centering\arraybackslash}p{#1}}
\definecolor{Fcolor}{rgb}{0, 0.5, 0.25}

\ifx\theorem\undefined
\newtheorem{theorem}{Theorem}
\newenvironment{theorem*}{\par\noindent{\bf Theorem\ }}{\hfill\\[2mm]}

\newtheorem{proposition}[theorem]{Proposition}

\newenvironment{corollary*}{\par\noindent{\bf Corollary\ }}{\hfill\\[2mm]}

\fi

\newcommand{\fs}{\sm{\mathrm{f}}{5}}
\newcommand{\ps}{\sm{\mathrm{p}}{5}}
\newcommand{\objV}{\mathrm{obj}}
\newcommand{\kV}{k}
\newcommand{\smallMinus}{\scalebox{1.2}[0.7]{-}}

\newcommand{\msh}{\!\!\;}

\DeclareMathAlphabet\EuRoman{U}{eur}{m}{n}
\SetMathAlphabet\EuRoman{bold}{U}{eur}{b}{n}
\newcommand{\euler}{\EuRoman}

\newcommand{\tr}{\mathrm{tr}}

\newcommand{\zerbf}{\mathbf{0}}

\newcommand{\abf}{\mathbf{a}}
\newcommand{\Abf}{\mathbf{A}}

\newcommand{\Abfc}{\boldsymbol{\euler{A}}}

\newcommand{\Bbf}{\mathbf{B}}

\newcommand{\Bbfc}{\boldsymbol{\euler{B}}}

\newcommand{\cbf}{\mathbf{c}}
\newcommand{\Ccal}{\mathcal{C}}
\newcommand{\Cbfc}{\boldsymbol{\euler{C}}}

\newcommand{\Dbfc}{\boldsymbol{\euler{D}}}

\newcommand{\Ebf}{\mathbf{E}}

\newcommand{\Fbf}{\mathbf{F}}

\newcommand{\Gbfc}{\euler{G}}

\newcommand{\hbfc}{\boldsymbol{\euler{h}}}
\newcommand{\hfc}{\euler{h}}
\newcommand{\Hbf}{\mathbf{H}}
\newcommand{\Hcal}{\mathcal{H}}
\newcommand{\Ibf}{\mathbf{I}}

\newcommand{\Kbf}{\mathbf{K}}

\newcommand{\Kbfc}{\boldsymbol{\euler{K}}}

\newcommand{\Lbf}{\mathbf{L}}

\newcommand{\Mbfc}{\boldsymbol{\euler{M}}}

\newcommand{\Nbb}{\mathbb{N}}

\newcommand{\Ncal}{\mathcal{N}}

\newcommand{\Pbfc}{\boldsymbol{\euler{P}}}

\newcommand{\Qbfc}{\boldsymbol{\euler{Q}}}

\newcommand{\Rbb}{\mathbb{R}}
\newcommand{\Rbfc}{\boldsymbol{\euler{R}}}

\newcommand{\sbf}{\mathbf{s}}

\newcommand{\Sbb}{\mathbb{S}}

\newcommand{\ubfc}{\boldsymbol{\euler{u}}}

\newcommand{\wbfc}{\boldsymbol{\euler{w}}}
\newcommand{\Wbfc}{\boldsymbol{\euler{W}}}

\newcommand{\xbf}{\mathbf{x}}
\newcommand{\xbfc}{\boldsymbol{\euler{x}}}
\newcommand{\Xbf}{\mathbf{X}}

\newcommand{\ybfc}{\boldsymbol{\euler{y}}}

\newcommand{\Ybfc}{\boldsymbol{\euler{Y}}}

\newcommand{\zbfc}{\boldsymbol{\euler{z}}}

\newcommand{\nxs}{n_{\msh\sm{x}{5}}}
\newcommand{\nys}{n_{\msh\sm{y}{5}}}
\newcommand{\nzs}{n_{\msh\sm{z}{5}}}
\newcommand{\nws}{n_{\msh\sm{w}{5}}}
\newcommand{\nus}{n_{\msh\sm{u}{5}}}

\newcommand{\sm}[2]{\scaleto{#1\mathstrut}{#2pt}}
\definecolor{pinkF}{rgb}{0.858, 0.188, 0.478}
\newif\ifcomment

\crefrangelabelformat{equation}{(#3#1#4)\,--\,(#5#2#6)}
\crefmultiformat{equation}{(#2#1#3)}{\,--\,(#2#1#3)}{#2#1#3}{#2#1#3}

\crefname{equation}{}{}
\crefname{figure}{Figure}{Figures}
\crefname{algorithm}{Algorithm}{}
\crefname{table}{Table}{Tables}
\crefname{lemma}{Lemma}{Lemmas}

\commentfalse
\IEEEoverridecommandlockouts                              
\title{\LARGE \bf Convex Relaxation of Bilinear Matrix Inequalities\\ Part II: Applications to Optimal Control Synthesis}
\author{Mohsen Kheirandishfard, Fariba Zohrizadeh, Muhammad Adil, and Ramtin Madani
\thanks{Mohsen Kheirandishfard and Fariba Zohrizadeh are with the Department of Computer Science and Engineering, The University of Texas at Arlington, Arlington, TX 76019, USA (emails:mohsen.kheirandishfard@uta.edu, fariba.zohrizadeh@uta.edu), Muhammad Adil and Ramtin Madani are with the Department of Electrical Engineering, The University of Texas at Arlington, Arlington, TX 76019, USA (emails:muhammad.adil@uta.edu, ramtin.madani@uta.edu). This work is in part supported by the NSF award 1809454 and a University of Texas System STARs award.}%
}

\begin{document}
\maketitle
\thispagestyle{empty}
\pagestyle{empty}
\begin{abstract}
The first part of this paper proposed a family of penalized convex relaxations for solving optimization problems with bilinear matrix inequality (BMI) constraints. In this part, we generalize our approach to a sequential scheme which starts from an arbitrary initial point (feasible or infeasible) and solves a sequence of penalized convex relaxations in order to find feasible and near-optimal solutions for BMI optimization problems. We evaluate the performance of the proposed method on the $\mathcal{H}_2$ and $\mathcal{H}_{\infty}$ optimal controller design problems with both centralized and decentralized structures. The experimental results based on a variety of benchmark control plants demonstrate the promising performance of the proposed approach in comparison with the existing methods. 
\end{abstract}

\section{Introduction}

The design of optimal controllers can be computationally challenging due to NP-hardness in general
\cite{tsitsiklis1985complexity,witsenhausen1968counterexample,papadimitriou1986intractable}.
This two part paper is motivated by the applications of bilinear matrix inequalities (BMIs) in optimal control. We study the class of optimization problems with BMI constraints and their applications in the design of optimal structured controllers. In Part I, we proposed a variety of convex relaxations and penalization methods for solving BMI optimization problems.
In this part, we investigate the $\Hcal_2$ and $\Hcal_{\infty}$ optimal structured controller design tasks by means of the proposed convex relaxations. 
Structured controllers has been extensively explored for various systems, including
for spatially distributed systems \cite{bamieh2002distributed,motee2008optimal,wu2017sparsity,d2003distributed,ehsan2018distributed}, localizable systems \cite{wang2014localized}, energy systems \cite{dorfler2014sparsity,schuler2014decentralized,zhang2017distributed}, optimal static distributed systems \cite{lin2011augmented,fardad2009optimal}, strongly connected systems \cite{lavaei2012decentralized}, as well as heterogeneous systems \cite{dullerud2004distributed}.

The problems of designing distributed state-feedback and output-feedback controllers for linear time-invariant systems have been studied by several papers \cite{bahavarnia2015sparse,lin2013design,dhingra2016method,scherer2000design,lavaei2008control,palacios2014recent,de2000convexifying,arastoo2015output}. 
The papers \cite{bamieh2002optimal,fan1994centralized,qi2004structured,voulgaris2001convex,hengster2015distributed} have considered special cases which make controller design problems computationally tractable. 
The paper \cite{rotkowitz2006characterization} introduces a condition regarded as quadratic invariance, which enables the transformation of optimal distributed controller design problems to convex optimization. This condition is further explored by other papers, including \cite{tanaka2014optimal,lessard2012optimal,lamperski2013output,rotkowitz2012nearest,shah2013cal,alavian2013q,matni2013dual,Matni2016regul}. 
Inspired by 
\cite{toker1995np,gross2011optimized,de2002extended,ikeda1996decentralized,sojoudi2010interconnection}, we cast a variety of controller design tasks as optimization subject to BMI constraints and employ a family of penalized convex relaxations to solve the resulting nonconvex optimization problems. 
Our work is related to the body of literature on convex relaxation of optimal controller design based on semidefinite programming in
\cite{fazelnia2017convex,fattahi2015transformation,lin2017convex,zhang2017distributed}, as well as sequential methods in \cite{wang2009time,zhai2001decentralized,sadabadi2013lmi,fardad2014design}.


We discuss the state-of-the-art semidefinite programming (SDP) and second-order cone programming (SOCP) relaxations. Additionally, we introduce a computationally efficient parabolic relaxation which solely relies on convex quadratic inequalities as opposed to conic constraints. Next, a family of penalty functions are introduced that can leverage any arbitrary initial point. The incorporation of these penalty terms into the objective of SDP, SOCP, and parabolic relaxations is guaranteed to produce feasible points for BMI optimization problems, as long as the initial point is sufficiently close the BMI feasible set. 

Built upon the theoretical results of Part I, in this part, we offer a sequential penalized relaxation which is able to find feasible and near-globally optimal solutions for BMI optimization. The proposed sequential method is applied to a variety of $\Hcal_2$ and $\Hcal_{\infty}$ static output-feedback controller design problems, and
its performance is tested on 
control plants from the COMPl\textsubscript{e}ib \cite{leibfritz2006compleib} library. 
Numerical experiments demonstrate the promising performance of the proposed method in comparison with the existing methods and software packages.

\subsection{Notation}
Throughout the paper, the scalars, vectors, and matrices are respectively shown by italic letters, lower-case bold letters, and upper-case bold letters. Symbols $\Rbb$, $\Rbb^{n}$, and $\Rbb^{n\times m}$ respectively denote the set of real scalars, real vectors of size $n$, and real matrices of size $n\times m$. The set of real $n\times n$ symmetric matrices is shown with $\Sbb_n$. For a given vector $\abf$ and matrix $\Abf$, symbols $a_i$ and $A_{ij}$ respectively indicate the $i^{th}$ element of $\abf$ and $(i,j)^{th}$ element of $\Abf$. For symmetric matrix $\Abf$, notations $\Abf\succeq 0$ and $\Abf\preceq 0$ show positive and negative semidefinite ($\Abf\succ 0$ and $\Abf\prec 0$ indicate positive and negative definite). For two arbitrary matrices $\Abf$ and $\Bbf$ of the same size, symbol $\langle\Abf,\Bbf\rangle=\tr\{\Abf^{\!\top}\Bbf\}$ shows the inner product between the matrices where $\tr\{.\}$ and $(.)^\top$ respectively denote the trace and transpose operators. Operator $\mathrm{diag}(.)$ gets a vector and forms a diagonal matrix with its input on the diagonal. For a symmetric matrix $\Bbf$ of size $n$, symbol $\Bbf(:)$ indicates a vector of size $\binom{n}{2}$ consists of all unique elements of $\Bbf$. Symbols $\Ibf$ and $\mathbf{0}$ respectively denote the identity matrix and zero matrix of appropriate dimensions and $\Ncal$ is a shorthand for the set $\{1,\dots,n\}$. 

\section{Problem Formulation}
In this part, we formulate the problem of structured static output-feedback controller design as an optimization problem with linear objective function and a bilinear matrix inequality (BMI) constraint, as follows: 
\begin{subequations} \label{eq:ODC_BMIOrig}
\begin{align}
&\underset{{\xbf}\in{\Rbb}^{n}}{\text{minimize}} && \cbf^{\!\top}\xbf \label{eq:ODC_BMIOrig_obj}\\
&\text{subject to}
&&p(\xbf,\xbf\xbf^{\!\top})\preceq 0,\label{eq:ODC_BMIOrig_con_01}
\end{align}
\end{subequations}
where $\cbf\in\Rbb^{n}$ is given, $p\!:\!\Rbb^{n}\!\times\! \Sbb_{n}\!\rightarrow\Sbb_{m}$ is a matrix-valued function defined as
\begin{equation}\label{eq:ODC_MatPen}
\begin{aligned}[b]
& p(\xbf,\Xbf)\triangleq{\Fbf}_{\sm{0}{6}}+\sum_{\sm{k\in\Ncal}{5.5}} {x}_{k}{\Kbf}_k+\sum_{\sm{i\in\Ncal}{5.5}}\sum_{\sm{j\in\Ncal}{5.5}}{X}_{ij}{\Lbf}_{ij},
\end{aligned}
\end{equation}
and ${\Fbf_{\sm{0}{6}}}$, $\{{\Kbf}_k\}_{\sm{k\in\Ncal}{5.5}}$, and $\{{\Lbf}_{ij}\}_{\sm{i,j\in\Ncal}{5.5}}$ are fixed $m\times m$ real symmetric matrices. Due to the presence of bilinear terms, the problem \cref{eq:ODC_BMIOrig_obj,eq:ODC_BMIOrig_con_01} is nonconvex and NP-hard in general. To address this issue, in Part I, we defined an auxiliary matrix variable $\Xbf\in\Sbb_n$ to account for $\xbf\xbf^{\top}$, and developed a computationally-tractable convex surrogate of the following form:
\begin{subequations} \label{eq:ODC_GenRelaxation_Pen}
	\begin{align}
	&\underset{{\xbf}\sm{\in}{6}{\Rbb}^{n}\!,{\Xbf}\sm{\in}{6}{\Sbb}_{n}}{\text{minimize}} && \cbf^{\!\top} \xbf+{\eta}\;(\mathrm{tr}\{\Xbf\}-2\;{\check{\xbf}}^{\!\top}{\xbf}+\check{\xbf}^{\!\top}{\check{\xbf}}) \label{eq:ODC_GenRelaxation_Pen_obj}\\
	&\text{\,subject to}
	&&p(\xbf,\Xbf)\;\preceq 0,\label{eq:ODC_GenRelaxation_Pen_con_01} \\
	&&& \Xbf-\xbf\xbf^{\!\top}\!\in\Ccal, \label{eq:ODC_GenRelaxation_Pen_con_02}
	\end{align}
\end{subequations}
where $\check{\xbf}$ is an arbitrary initial guess, $\Ccal$ is a proper convex cone, and $\eta\!>\!0$ is a regularization parameter which offers a trade-off between the linear objective function and the penalty term. Observe that the BMI constraint \cref{eq:ODC_BMIOrig_con_01} is transformed to the linear matrix inequality (LMI) constraint \cref{eq:ODC_GenRelaxation_Pen_con_01} with respect to $\xbf$ and $\Xbf$. Additionally, the relation $\Xbf=\xbf\xbf^{\!\top}$ is relaxed to the convex constraint \cref{eq:ODC_GenRelaxation_Pen_con_02}. In Part I, the following three choices for the convex cone $\mathcal{C}$ are discussed
\begin{subequations}
	\begin{align}
	&\hspace{-0.2cm}\Ccal_{\sm{\mathrm{1}}{5}}\!\!\!&&\hspace{-0.26cm}=\!\{\Hbf\!\in\!\Sbb_{n}\!\mid\!\Hbf\!\!\!\!\!&&&&\hspace{-0.5cm}\succeq 0\},\\[2.5pt]
	&\hspace{-0.2cm}\Ccal_{\sm{\mathrm{2}}{5}}\!\!\!&&\hspace{-0.26cm}=\!\{\Hbf\!\in\!\Sbb_{n}\!\mid\! H_{ii}\!\!\!\!\!&&&&\hspace{-0.5cm}\geq 0,H_{ii}H_{jj}\geq H_{ij}^2,\sm{\forall i,\!j\!\in\!\Ncal}{8}\}, \\[2.5pt]
	&\hspace{-0.2cm}\Ccal_{\sm{\mathrm{3}}{5}}\!\!\!&&\hspace{-0.26cm}=\!\{\Hbf\!\in\!\Sbb_{n}\!\mid\! H_{ii}\!\!\!\!\!&&&&\hspace{-0.5cm}\geq 0,H_{ii}\!+\!H_{jj}\!\geq\! 2\left|H_{ij}\right|\!,\sm{\forall i,\!j\!\in\!\Ncal}{8}\},
    \end{align}
\end{subequations}
leading to the semidefinite programming (SDP), second-order cone programming (SOCP), and parabolic relaxations, respectively. The resulting convex relaxations are efficiently solvable up to any desired accuracy using the standard numerical algorithms. 

In the next section, we propose a sequential framework which solves penalized convex relaxations of form \cref{eq:ODC_GenRelaxation_Pen_obj,eq:ODC_GenRelaxation_Pen_con_02} to obtain feasible and near-globally optimal points for the original BMI problem \cref{eq:ODC_BMIOrig_obj,eq:ODC_BMIOrig_con_01}.

\section{Sequential Penalized Relaxation}
In Part I, it is proven that the proposed penalized convex relaxation is guaranteed to preserve the feasibility of any Mangasarian-Fromovitz regular initial point. Moreover, we proved that infeasible initial points that are sufficiently close to the BMI feasible set enjoy a similar property. In light of these theoretical guarantees, we propose an algorithm which can start from an arbitrary initial point and proceed sequentially until a satisfactory solution of the BMI problem is obtained. Once feasibility is attained, it is preserved and the objective value is improved in the subsequent round of algorithm. This procedure is detailed in \cref{al:ODC_alg_01}. As it is shown, the algorithm is stopped if the number of rounds exceeds $\mathrm{maxRound}$ or if the improvement between two consecutive rounds is less than $\mathrm{progThresh}$. We observed that the Nesterov's acceleration \cite{nesterov1983method} (line 8 of \cref{al:ODC_alg_01}) can considerably improve the convergence behavior of the proposed algorithm.
\begin{algorithm}
	\caption{Sequential Penalized Relaxation}
	\label{SP-SDP}
	\begin{algorithmic}[1]
		\Require {$\check{\xbf}\msh\in\msh\Rbb^{n}$, $\eta\msh>\msh0,\mathrm{progThresh}\msh\in\msh\Rbb$, $\mathrm{maxRound}\msh\in\msh\Nbb$} 
		\Ensure {$\accentset{\ast}{\xbf}$}
        \State $\xbf_0\phantom{k}\!\!\!\!\leftarrow\check{\xbf}$
        \State $k\phantom{\xbf_0}\!\!\!\!\leftarrow 1$
		\Repeat
			\State $\xbf_k\leftarrow$ Solve penalized relaxation \cref{eq:ODC_GenRelaxation_Pen_obj,eq:ODC_GenRelaxation_Pen_con_01,eq:ODC_GenRelaxation_Pen_con_02}
		\If {$\frac{\lvert\cbf^{\!\top}(\xbf_k-\xbf_{k-1})\rvert}{\phantom{\frac{1}{2}}\!\!\!\lvert\cbf^{\!\top}\xbf_{k-1}\rvert}\!\times\! 100\leq \mathrm{progThresh}$}\vspace{0.04cm}
		\State \textbf{break}
		\EndIf
        \State $\check{\xbf} \leftarrow\xbf_k+\frac{k-1}{k+2}(\xbf_k-\xbf_{k-1})$
        \State $k\leftarrow k+1$
		\Until $k\leq\mathrm{maxRound}$
		\State $\accentset{\ast}{\xbf}\leftarrow\xbf_{k-1}$
	\end{algorithmic}\label{al:ODC_alg_01}
\end{algorithm}

The next section is concerned with the formulation of the $\mathcal{H}_2$ and $\mathcal{H}_{\infty}$ optimal structured control synthesis problems in the form the optimization \cref{eq:ODC_BMIOrig_obj,eq:ODC_BMIOrig_con_01}.

\section{Optimal Structured Controller Synthesis}
Consider a linear-time invariant control plant $\Gbfc$ with the following dynamics:

\begin{equation}  \label{eq:ODC_PlantOl}
\begin{aligned}
\Gbfc:\begin{cases}
&\hspace{-0.2cm}\dot{\xbfc}=\\
&\hspace{-0.2cm}\zbfc=\\
&\hspace{-0.2cm}\ybfc=
\end{cases}
\end{aligned}
\begin{aligned}
&\Abfc&&\hspace{-0.3cm}\xbfc+\Bbfc_{\sm{1}{6}}\\
&\Cbfc_{\sm{1}{6}}&&\hspace{-0.3cm}\xbfc+\Dbfc_{\sm{11}{6}}\\
&\Cbfc&&\hspace{-0.3cm}\xbfc+\Dbfc_{\sm{21}{6}}
\end{aligned}
\begin{aligned}
&\hspace{-0.05cm}\wbfc+\Bbfc&&\hspace{-0.33cm}\ubfc\\
&\hspace{-0.05cm}\wbfc+\Dbfc_{\sm{12}{6}}&&\hspace{-0.33cm}\ubfc\\
&\hspace{-0.05cm}\wbfc\phantom{+\Dbfc_{\sm{12}{6}}}&&\hspace{-0.33cm}
\end{aligned}
\end{equation}
where $\xbfc\in\Rbb^{\nxs}$ is the vector of states, $\wbfc\in\Rbb^{\nws}$ is the system inputs, $\ubfc\in\Rbb^{\nus}$ denotes the control commands, $\zbfc\in\Rbb^{\nzs}$ is the response outputs, and $\ybfc\in\Rbb^{\nys}$ indicate the vector of sensor measurements. The matrices $\Abfc$, $\Bbfc_{\sm{1}{6}}$, $\Bbfc$, $\Cbfc_{\sm{1}{6}}$, $\Cbfc$, $\Dbfc_{\sm{11}{6}}$, $\Dbfc_{\sm{12}{6}}$, $\Cbfc$, $\Dbfc_{\sm{21}{6}}$ are all fixed and of appropriate dimensions. In what follows, we aim to design a structured static output-feedback controller for this plant. Define $\Kbfc:\Rbb^l\rightarrow\Rbb^{\nus\times\nys}$ as the controller with the following structure:
\begin{equation}
\begin{aligned}
\Kbfc(\msh\hbfc\msh)\triangleq \sum_{i=1}^l \hfc_i \Ebf_i,
\end{aligned}
\end{equation}
where $\hbfc\in\Rbb^l$ represents nonzero elements of the controller and $\{\Ebf_i\}_{i=1}^l\in\{0,1\}^{\nus\times\nys}$ are fixed binary matrices. Given a vector $\ybfc$ of all measurements as the controller input, the controller outputs a command vector is given by the equation $\ubfc=\Kbfc(\msh\hbfc\msh)\ybfc$. In order to obtain the optimal controller, we first derive the dynamic equations that represent the closed-loop system $\Gbfc$ as:
\begin{equation}\label{eq:ODC_PlantCl}
\Gbfc_{\sm{\mathrm{cl}}{6}}:
\begin{aligned}
\begin{cases}
& \\
& 
\end{cases}
\end{aligned}
\begin{aligned}
&\hspace{-0.4cm}\dot{\xbfc}&&\hspace{-0.4cm}=\Abfc_{\sm{\mathrm{cl}}{6}}(\msh\hbfc\msh)\xbfc+\Bbfc_{\sm{\mathrm{cl}}{6}}(\msh\hbfc\msh)&&&&\hspace{-0.7cm}\wbfc \\
&\hspace{-0.4cm}\zbfc&&\hspace{-0.4cm}=\Cbfc_{\sm{\mathrm{cl}}{6}}(\msh\hbfc\msh)\xbfc+\Dbfc_{\sm{\mathrm{cl}}{6}}(\msh\hbfc\msh)&&&&\hspace{-0.7cm}\wbfc
\end{aligned}
\end{equation}
where the matrix functions 
$\Abfc_{\sm{\mathrm{cl}}{6}}\!:\!\Rbb^l\!\!\rightarrow\!\Rbb^{\nxs\!\times\nxs}$, $\Bbfc_{\sm{\mathrm{cl}}{6}}\!:\!\Rbb^l\!\!\rightarrow\!\Rbb^{\nxs\!\times\nws}$, $\Cbfc_{\sm{\mathrm{cl}}{6}}\!:\!\Rbb^l\!\!\rightarrow\!\Rbb^{\nzs\!\times\nxs}$, and $\Dbfc_{\sm{\mathrm{cl}}{6}}\!:\!\Rbb^l\!\!\rightarrow\!\Rbb^{\nzs\!\times\nws}$ are defined as
\begin{subequations} \label{eq:ODC_ClMa}
\begin{align}
& \phantom{\Bbfc_{\sm{\mathrm{cl}}{6}}\Cbfc_{\sm{\mathrm{cl}}{6}}\Dbfc_{\sm{\mathrm{cl}}{6}}}\Abfc_{\sm{\mathrm{cl}}{6}}(\msh\hbfc\msh)\triangleq\Abfc_{\phantom{\sm{11}{6}}}\phantom{\Bbfc_{\sm{1}{6}\phantom{\sm{1}{6}}}\Cbfc_{\sm{1}{6}\phantom{\sm{1}{6}}}\Dbfc_{\sm{11}{6}}}\hspace{-1.65cm}+\Bbfc_{\phantom{\sm{12}{6}}}\Kbfc(\msh\hbfc\msh)\Cbfc_{\phantom{\sm{21}{6}}},\label{eq:ODC_ClLoMatA} \\[1pt]
& \phantom{\Abfc_{\sm{\mathrm{cl}}{6}}\Cbfc_{\sm{\mathrm{cl}}{6}}\Dbfc_{\sm{\mathrm{cl}}{6}}}\Bbfc_{\sm{\mathrm{cl}}{6}}(\msh\hbfc\msh)\triangleq\Bbfc_{\sm{1}{6}\phantom{\sm{1}{6}}}\phantom{\Abfc_{\phantom{\sm{11}{6}}}\Cbfc_{\sm{1}{6}\phantom{\sm{1}{6}}}\Dbfc_{\sm{11}{6}}}\hspace{-1.65cm}+\Bbfc_{\phantom{\sm{12}{6}}}\Kbfc(\msh\hbfc\msh)\Dbfc_{\sm{21}{6}}, \label{eq:ODC_ClLoMatB} \\[1pt]
& \phantom{\Bbfc_{\sm{\mathrm{cl}}{6}}\Abfc_{\sm{\mathrm{cl}}{6}}\Dbfc_{\sm{\mathrm{cl}}{6}}}\Cbfc_{\sm{\mathrm{cl}}{6}}(\msh\hbfc\msh)\triangleq\Cbfc_{\sm{1}{6}\phantom{\sm{1}{6}}}\phantom{\Abfc_{\phantom{\sm{11}{6}}}\Bbfc_{\sm{1}{6}\phantom{\sm{1}{6}}}\Dbfc_{\sm{11}{6}}}\hspace{-1.65cm}+\Dbfc_{\sm{12}{6}}\Kbfc(\msh\hbfc\msh)\Cbfc_{\phantom{\sm{21}{6}}},\label{eq:ODC_ClLoMatC} \\[1pt]
& \phantom{\Bbfc_{\sm{\mathrm{cl}}{6}}\Cbfc_{\sm{\mathrm{cl}}{6}}\Abfc_{\sm{\mathrm{cl}}{6}}}\Dbfc_{\sm{\mathrm{cl}}{6}}(\msh\hbfc\msh)\triangleq\Dbfc_{\sm{11}{6}}\phantom{\Abfc_{\phantom{\sm{11}{6}}}\Bbfc_{\sm{1}{6}\phantom{\sm{1}{6}}}\Cbfc_{\sm{1}{6}\phantom{\sm{1}{6}}}}\hspace{-1.65cm}+\Dbfc_{\sm{12}{6}}\Kbfc(\msh\hbfc\msh)\Dbfc_{\sm{21}{6}}.\hspace{1cm}\label{eq:ODC_ClLoMatD}
\end{align}
\end{subequations}
In what follows, we cast the $\Hcal_2$- and $\Hcal_{\infty}$-norm optimal controller design problems for the plant $\Gbfc$ as BMI problems of form \cref{eq:ODC_BMIOrig_obj,eq:ODC_BMIOrig_con_01}.

\subsection{$\Hcal_2$ Optimal Control}
The $\Hcal_2$-norm of a control system is described as the average energy of the output signal, given white noise as the input. For the control plant $\Gbfc$, the $\Hcal_2$ optimal controller design problem aims to find a vector $\hbfc$ such that the structured controller $\Kbfc(\msh\hbfc\msh)$ stabilizes the plant (i.e. all the eigenvalues of $\Abfc_{\sm{\mathrm{cl}}{6}}(\msh\hbfc\msh)$ have negative real part) and minimizes the $\Hcal_2$ norm of $\Gbfc_{\sm{\mathrm{cl}}{6}}$. With no loss of generality, we assume that $\Dbfc_{\sm{11}{6}}\!=\!\mathbf{0}$, $\Dbfc_{\sm{21}{6}}\!=\!\mathbf{0}$, and that there exists a stabilizing controller gain $\Kbfc(\msh\hbfc\msh)$ for $\Gbfc$. Hence, the $\Hcal_2$ norm of the closed-loop plant is given as
\begin{equation}
\begin{aligned}
\lVert\Gbfc_{\sm{\mathrm{cl}}{6}}{\rVert}_{\Hcal_2}=\mathrm{tr}\{\Cbfc_{\sm{\mathrm{cl}}{6}}(\msh\hbfc\msh)\Pbfc\Cbfc_{\sm{\mathrm{cl}}{6}}(\msh\hbfc\msh)^{\!\top}\},
\end{aligned}
\end{equation}
where matrix $\Pbfc\succ 0$ is the solution of the following Lyapunov equation:
\begin{equation}\label{eq:ODC_Lyap_011}
\begin{aligned}
\Abfc_{\sm{\mathrm{cl}}{6}}(\msh\hbfc\msh)\Pbfc+\Pbfc\Abfc_{\sm{\mathrm{cl}}{6}}(\msh\hbfc\msh)^{\!\top}+\Bbfc_{\sm{1}{6}}\Bbfc_{\sm{1}{6}}^{\!\top}=0.
\end{aligned}
\end{equation}
It is well-known that, $\Pbfc$ can be obtained by solving the following relaxed matrix inequality \cite{boyd2004convex}:
\begin{equation}\label{eq:BMI_Lyap}
\begin{aligned}[b]
\Abfc_{\sm{\mathrm{cl}}{6}}(\msh\hbfc\msh)\Pbfc+\Pbfc\Abfc_{\sm{\mathrm{cl}}{6}}(\msh\hbfc\msh)^{\!\top}+\Bbfc_{\sm{1}{6}}\Bbfc_{\sm{1}{6}}^{\!\top}\preceq 0.
\end{aligned}
\end{equation}
Therefore, the following optimization problem minimizes the $\Hcal_2$ norm, subject to a stabilizing controller with a desired zero-nonzero pattern:
\begin{subequations} \label{eq:ODC_Htwo}
\begin{align}
&\underset{\begin{subarray}{c} {\Pbfc}\in{\Sbb}_{\nxs}\msh\msh,\msh\Wbfc\in\Sbb_{\nzs}\\{\hbfc}\in{\Rbb}^{l}\end{subarray}}{\text{minimize}} && \langle\Wbfc,\Ibf\rangle \label{eq:ODC_Htwo_obj}\\
&\text{~~subject to}
&& f_{\sm{\mathrm{LMI}}{4}}(\Pbfc,\Wbfc)+f_{\sm{\mathrm{BMI}}{4}}(\Pbfc,\hbfc)\preceq 0,\label{eq:ODC_Htwo_con_01}
\end{align}
\end{subequations}
where the matrix functions $f_{\sm{\mathrm{LMI}}{4}}\!:\!\Sbb_{\nxs}\!\times\!\Sbb_{\nzs}\!\rightarrow\!\Sbb_{2\nxs+\nzs}$ and $f_{\sm{\mathrm{BMI}}{4}}\!:\! \Sbb_{\nxs}\!\times\!\Rbb^{l}\!\rightarrow\!\Sbb_{2\nxs+\nzs}$ are defined as
\begin{subequations}
\begin{align}
&f_{\sm{\mathrm{LMI}}{4}}(\Pbfc,\Wbfc)\triangleq\nonumber\\
					&\begin{bmatrix} \Abfc\Pbfc+\Pbfc\Abfc^{\!\top}\!+\!{\Bbfc_{\sm{1}{6}}}{\Bbfc_{\sm{1}{6}}}^{\!\!\!\top}\phantom{a\,} & \zerbf & \zerbf \\
                    				\ast &  -\Wbfc & {\Cbfc_{\sm{1}{6}}}\Pbfc \\
                                    \ast & \ast & -\Pbfc \end{bmatrix},\label{eq:ODC_Htwo_f_LMI}\\
&f_{\sm{\mathrm{BMI}}{4}}(\Pbfc,\hbfc)\triangleq\nonumber\\
					&\begin{bmatrix} \Bbfc\Kbfc(\msh\hbfc\msh)\Cbfc\Pbfc\!+\!{(\Bbfc\Kbfc(\msh\hbfc\msh)\Cbfc\Pbfc)}^{\!\top} & \zerbf & \zerbf \\
                    				\ast &  \zerbf & \Dbfc_{\sm{12}{6}}\Kbfc(\msh\hbfc\msh)\Cbfc\Pbfc \\
                                    \ast & \ast & \zerbf \end{bmatrix}, \label{eq:ODC_Htwo_f_BMI}      
\end{align}
\end{subequations}
and $\ast$ account for the symmetric elements of the matrices.
\vspace{2mm}
\begin{proposition}
Assume that $\Bbfc_{\sm{1}{6}}\Bbfc_{\sm{1}{6}}^{\!\top}\succ 0$, $\Dbfc_{\sm{11}{6}}\!=\!\mathbf{0}$, and $\Dbfc_{\sm{21}{6}}\!=\!\mathbf{0}$. If $(\accentset{\ast}{\Pbfc},\accentset{\ast}{\Wbfc},\accentset{\ast}{\hbfc})$ is an optimal solution of the problem \cref{eq:ODC_Htwo_obj,eq:ODC_Htwo_con_01}, then $\Kbfc(\msh\accentset{\ast}{\hbfc}\msh)$ is the optimal $\Hcal_2$ static output-feedback controller gain for the plant $\Gbfc$.
\end{proposition}

\vspace{2mm}
\begin{proof}
It can be easily verified from the BMI constraint \cref{eq:ODC_Htwo_con_01}, that $\accentset{\ast}{\Pbfc}$ is positive-definite, and satisfies the Lyapunov inequality \cref{eq:BMI_Lyap}, which certifies that $\Kbfc(\msh\accentset{\ast}{\hbfc}\msh)$ is a stabilizing controller. On the other hand, we have $\accentset{\ast}{\Wbfc}=\Cbfc_{\sm{\mathrm{cl}}{6}}(\msh\accentset{\ast}{\hbfc}\msh)\accentset{\ast}{\Pbfc}\Cbfc_{\sm{\mathrm{cl}}{6}}(\msh\accentset{\ast}{\hbfc}\msh)^{\!\top}$ which means that the closed-loop gain $\langle\accentset{\ast}{\Wbfc},\Ibf\rangle=\lVert\Gbfc_{\sm{\mathrm{cl}}{6}}{\rVert}_{\Hcal_2}$ is minimized subject to the stability condition.
\end{proof} 
\vspace{2mm}

The constraint \cref{eq:ODC_Htwo_con_01} is a BMI due to the presence of matrix product $\Kbfc(\msh\hbfc\msh)\Cbfc\Pbfc$. Hence, it can be easily observed that the controller design problem \cref{eq:ODC_Htwo_obj,eq:ODC_Htwo_con_01} can be cast in the form of \cref{eq:ODC_BMIOrig_obj,eq:ODC_BMIOrig_con_01}. 
To this end, we stack all of the variables into a vector
\begin{equation}\label{eq:ODC_Htwo_xdef}
\begin{aligned}
\tilde{\xbf}\triangleq\mathrm{diag}(\tilde{\sbf})[\Wbfc(:)^{\!\top},\Pbfc(:)^{\!\top},\hbfc^{\!\top}]^{\top}\in\Rbb^{\tilde{n}},
\end{aligned}
\end{equation}
where $\tilde{n}\!=\!\binom{\nxs}{2}+\binom{\nzs}{2}+l$ and $\tilde{\sbf}\in\Rbb^{\tilde{n}}$ 
is fixed. Since the performance of the proposed sequential algorithm depends on the choice of coordinates, we consider the following values for the elements of the scale vector $\tilde{\sbf}$:
\begin{itemize}
\item ${\tilde{s}}_i=1$ if $\exists j\in\Ncal$, such that $\tilde{x}_i\tilde{x}_j$ appears in the BMI constraint \cref{eq:ODC_BMIOrig_con_01}.
\item ${\tilde{s}}_i=\mathrm{min}(0.5\eta,0.01)$ if $\nexists j\in\Ncal$, such that $\tilde{x}_i\tilde{x}_j$ appears in the BMI constraint \cref{eq:ODC_BMIOrig_con_01}.
\end{itemize}

\subsection{$\Hcal_{\infty}$ Optimal Control}
For a general control system, the $\Hcal_{\infty}$-norm is defined as the maximal amplification (system gain) from the input signal to the output. The problem of $\Hcal_{\infty}$ controller design for the linear system $\Gbfc$ aims at finding a vector $\hbfc$ such that the controller gain $\Kbfc(\msh\hbfc\msh)$ stabilizes the plant (i.e. all the eigenvalues of $\Abfc_{\sm{\mathrm{cl}}{6}}(\msh\hbfc\msh)$ have negative real part) and minimizes the $\Hcal_{\infty}$ norm of the closed-loop system. Assume that $\Dbfc_{\sm{21}{5}}=\zerbf$, $\gamma>0$, and that there exists a stabilizing controller $\Kbfc(\msh\hbfc\msh)$. Then, $\lVert\Gbfc_{\sm{\mathrm{cl}}{6}}{\rVert}_{\Hcal_{\infty}}<\gamma$ holds true if there exist a unique matrix $\Ybfc\succeq 0$ and controller $\Kbfc(\msh\hbfc\msh)$ that satisfy the following algebraic Riccati equation \cite{leibfritz2006compleib}:
\begin{equation}
\begin{aligned}
&\Abfc_{\sm{\mathrm{cl}}{6}}(\msh\hbfc\msh)\Ybfc+\Ybfc\Abfc_{\sm{\mathrm{cl}}{6}}(\msh\hbfc\msh)^{\!\top}+\gamma^{\smallMinus 1}\Bbfc_{\sm{\mathrm{cl}}{6}}(\msh\hbfc\msh)\Bbfc_{\sm{\mathrm{cl}}{6}}(\msh\hbfc\msh)^{\!\top}+\\
&{\gamma}^{\smallMinus 1}\Mbfc(\Ybfc,\!\!\:\hbfc,\!\!\:\gamma{)}^{\!\top}\!\Rbfc(\hbfc,\gamma)^{\smallMinus 1}\Mbfc(\Ybfc,\!\!\:\hbfc,\!\!\:\gamma)=0,
\end{aligned}
\end{equation}
where the matrix functions $\Mbfc : \Sbb_{\nxs}\!\!\times\!\Rbb^l\!\times\!\Rbb\!\rightarrow\!\Rbb^{\nzs\times\nxs}$ and $\Rbfc : \Rbb^l\!\times\!\Rbb\!\rightarrow\!\Sbb_{\nzs}$ are defined as
\begin{subequations}
\begin{align}
&\Rbfc(\hbfc,\gamma) \hspace{-0.75cm}&& \triangleq \Ibf-\gamma^{\smallMinus 2}\Dbfc_{\sm{\mathrm{cl}}{6}}(\msh\hbfc\msh)\Dbfc_{\sm{\mathrm{cl}}{6}}(\msh\hbfc\msh)^{\!\top},\\
&\Mbfc(\Ybfc,\!\!\:\hbfc,\!\!\:\gamma) \hspace{-0.75cm}&& \triangleq \Cbfc_{\sm{\mathrm{cl}}{6}}(\msh\hbfc\msh)\Ybfc+\gamma^{\smallMinus 1}\Dbfc_{\sm{\mathrm{cl}}{6}}(\msh\hbfc\msh)\Bbfc_{\sm{\mathrm{cl}}{6}}(\msh\hbfc\msh)^{\!\top}.
\end{align}
\end{subequations}
The existence of such solution is guaranteed if there exist $\Qbfc\succ\Ybfc\succeq 0$ and $\Kbfc(\msh\hbfc\msh)$ that satisfy 
\begin{equation}\label{eq:ODC_HinfCondO}
\begin{aligned}
&\Abfc_{\sm{\mathrm{cl}}{6}}(\msh\hbfc\msh)\Qbfc+\Qbfc\Abfc_{\sm{\mathrm{cl}}{6}}(\msh\hbfc\msh)^{\!\top}+\gamma^{\smallMinus 1}\Bbfc_{\sm{\mathrm{cl}}{6}}(\msh\hbfc\msh)\Bbfc_{\sm{\mathrm{cl}}{6}}(\msh\hbfc\msh)^{\!\top}+\\
&{\gamma}^{\smallMinus 1}\Mbfc(\Qbfc,\!\!\:\hbfc,\!\!\:\gamma{)}^{\!\top}\!\Rbfc(\hbfc,\gamma)^{\smallMinus 1}\Mbfc(\Qbfc,\!\!\:\hbfc,\!\!\:\gamma)\prec 0.
\end{aligned}
\end{equation}
Using Schur complement, the $\Hcal_{\infty}$ control design problem can be cast as the following optimization problem:
\begin{subequations} \label{eq:ODC_Hinf}
\begin{align}
&\underset{\begin{subarray}{c} {\Qbfc}\in{\Sbb}_{\nxs},\gamma\in\Rbb\\ {\hbfc}\in{\Rbb}^{l}\end{subarray}}{\text{minimize}} && \gamma \label{eq:ODC_Hinf_obj}\\
&\text{~subject to}
&&g_{\sm{\mathrm{LMI}}{4}}(\Qbfc,\gamma)+g_{\sm{\mathrm{BMI}}{4}}(\Qbfc,\hbfc)\preceq 0, \label{eq:ODC_Hinf_con_01}
\end{align}
\end{subequations}
where the matrix functions $g_{\sm{\mathrm{LMI}}{4}}:\Sbb_{\nxs}\!\!\times\!\Rbb\rightarrow\Sbb_{2\nxs+\nws+\nzs}$ and $g_{\sm{\mathrm{BMI}}{4}}:\Sbb_{\nxs}\!\!\times\!\Rbb^l\!\rightarrow\!\Sbb_{2\nxs+\nws+\nzs}$ are defined as
\begin{subequations}
\begin{align}
&\hspace{-0.1cm}g_{\sm{\mathrm{LMI}}{4}}(\Qbfc,\gamma)\!\triangleq\nonumber\\
						&\hspace{-0.3cm}\begin{bmatrix} -\Qbfc & \zerbf 									& \zerbf						 & \zerbf \\
                        				\ast   & \Abfc\Qbfc\!+\!\Qbfc\Abfc^{\!\!\top} & (\Cbfc_{\sm{1}{5}}\Qbfc)^{\!\top} & \Bbfc_{\sm{1}{6}}\phantom{\;} \\
                                        \ast   & \ast 								& -\gamma\Ibf 			 & \Dbfc_{\sm{11}{5}} \\
                                        \ast & \ast & \ast & -\gamma\Ibf \end{bmatrix}\!\!, \!\!\!\label{eq:ODC_Hinf_g_LMI}\\
&\hspace{-0.1cm}g_{\sm{\mathrm{BMI}}{4}}(\Qbfc,\hbfc)\!\triangleq\nonumber\\
						&\hspace{-0.3cm}\begin{bmatrix} \zerbf & \!\!\zerbf & \!\!\!\zerbf & \!\!\!\zerbf \\ 
                        \ast & \!\!\Bbfc\Kbfc(\msh\hbfc\msh)\Cbfc\Qbfc\!+\!(\Bbfc\Kbfc(\msh\hbfc\msh)\Cbfc\Qbfc)^{\!\top} & \!\!\!(\Dbfc_{\sm{12}{5}}\Kbfc(\msh\hbfc\msh)\Cbfc\Qbfc)^{\!\top} & \!\!\!\zerbf \\ 
                        \ast & \!\!\ast & \!\!\!\zerbf & \!\!\!\zerbf \\ 
                        \ast & \!\!\ast & \!\!\!\ast & \!\!\!\zerbf \end{bmatrix}\!\!. \!\!\!\label{eq:ODC_Hinf_g_BMI}
\end{align} 
\end{subequations}
\vspace{2mm}

\begin{proposition}
Assume that $\Dbfc_{\sm{11}{5}}\neq\zerbf$ and $\Bbfc_{\sm{1}{6}}\Bbfc_{\sm{1}{6}}^{\!\top}\succ 0$. If $(\accentset{\ast}{\Qbfc},\accentset{\ast}{\hbfc},\accentset{\ast}{\gamma})$ is an optimal solution of problem \cref{eq:ODC_Hinf_obj,eq:ODC_Hinf_con_01}, then $\Kbfc(\accentset{\ast}{\hbfc})$ is the optimal $\Hcal_{\infty}$ static output-feedback controller gain for plant $\Gbfc$.
\end{proposition}
\vspace{2mm}

\begin{proof}
From the BMI constraint \cref{eq:ODC_Hinf_con_01}, it can be easily verified that the assumption $\Dbfc_{\sm{11}{5}}\neq\zerbf$ concludes $\accentset{\ast}{\gamma}>0$. Moreover, the matrix $\accentset{\ast}{\Qbfc}$ is positive-definite and satisfies the inequality \cref{eq:ODC_HinfCondO}, which certifies that $\Kbfc(\accentset{\ast}{\hbfc})$ is a stabilizing controller. On the other hand, $\accentset{\ast}{\gamma}$ is minimized by the optimization problem \cref{eq:ODC_Hinf_obj,eq:ODC_Hinf_con_01}.
\end{proof}
\vspace{2mm}

Observe that \cref{eq:ODC_Hinf_con_01} is a BMI constraint as well, because of the matrix product $\Qbfc\Bbfc\Kbfc(\msh\hbfc\msh)$. Hence, we can cast \cref{eq:ODC_Hinf_obj,eq:ODC_Hinf_con_01} as an optimization problem of form \cref{eq:ODC_BMIOrig_obj,eq:ODC_BMIOrig_con_01}, with respect to the vector
\begin{equation}\label{eq:ODC_Hinf_xdef}
\begin{aligned}
\bar{\xbf}\triangleq \mathrm{diag}(\bar{\sbf})[\Qbfc(:)^{\!\top},\hbfc^{\!\top},\gamma]^{\top}\in\Rbb^{\bar{n}},
\end{aligned}
\end{equation}
where $\bar{n}\!=\!\binom{\nxs}{2}\!+\!l\!+\!1$ and the scale vector $\bar{\sbf}\in\Rbb^{\bar{n}}$ is created in a similar way as the $\Hcal_2$ case.

Next, we use \cref{al:ODC_alg_01} to find feasible and near-optimal solutions for the $\Hcal_2$ and $\Hcal_{\infty}$ optimal structured control design problems.

\section{Experimental Results}
In this section, the effectiveness of \cref{al:ODC_alg_01} is tested through extensive experiments on benchmark control plants from COMPl\textsubscript{e}ib \cite{leibfritz2006compleib}. The test cases cover a variety of applications, such as aircraft models (AC), academic test problems (NN), and decentralized interconnected systems (DIS), etc. We investigate the $\Hcal_2$ and $\Hcal_{\infty}$ optimal controller design problems for plants that are inherently static output-feedback stabilizable. 

We use the HIFOO  \cite{arzelier2011h2,burke2006hifoo} and the PENBMI \cite{kocvara2005penbmi} packages as competing solvers. The HIFOO is a publicly available MATLAB package which is based on a two-stage method for solving fixed order $\Hcal_2$ and $\Hcal_{\infty}$ output-feedback controller design problems. The first stage relies on the standard Broyden-Fletcher-Goldfarb-Shanno (BFGS) algorithm, and the second stage is based on random gradient sampling. The PENBMI package is a commercial local optimization solver, that is able to handle general BMI constrained problems with quadratic objectives. In our experiments, we have initialized all of the solvers with zero input. Other parameters of HIFOO and PENBMI are set to their default values. The experiments are all performed in MATLAB 2017a on a desktop computer with a 4-core 3.6GHz CPU and 32GB RAM. MOSEK v7 \cite{mosek2015mosek} is used through CVX to solve the resulting convex programs.

The reminder of this section offers detailed discussion of our experiments on centralized and fully decentralized controller design problems. 
\begin{table*}[!t] 
\vspace{0.09cm}
\footnotesize
	\centering
	\caption{Results of centralized $\Hcal_2$ controller design for COMPl\textsubscript{e}ib models.}
    \scalebox{0.95}{
	\begin{tabular}{ |@{\;}c@{\;}|@{\;}c@{\;}
    |@{\;}|
 @{\!\!\:\,}c@{\;}|@{\;}c@{\;}|@{\;}c@{\;}|@{\;}c@{\;}|@{\;}c@{\;}|@{\;}c@{\;}|@{\;}|@{\;}c@{\;}|@{\;}c@{\;}|@{\;}c@{\;}|@{\;}c@{\;}|@{\;}c@{\;}|@{\;}c@{\;}
 |@{\;}|
 @{\;}c@{\;}|@{\;}c@{\;}|@{\;}c@{\;}|@{\;}c@{\;}|@{\;}c@{\;}|@{\;}c@{\;}
 |@{\;}|
 @{\;}c@{\;}|@{\;}c@{\;}|}
		\hline
		\multirow{2}{*}{Name} & \multirow{2}{*}{$\lVert\Gbfc_{\sm{\mathrm{ol}}{6}}{\rVert}_{2}$} & \multicolumn{6}{c|@{\;}|@{\;}}{SDP} & \multicolumn{6}{c|@{\;}|@{\;}}{SOCP} & \multicolumn{6}{c|@{\;}|@{\;}}{Parabolic} & \multicolumn{2}{c|}{Competitors}\\
		\cline{3-22}
 		& &  $\eta$ & t & ${\kV}_{\fs}$ & ${\objV}_{\fs}$ & ${\kV}_{\ps}$ & ${\objV}_{\ps}$ & $\eta$ & t & ${\kV}_{\fs}$ & ${\objV}_{\fs}$ & ${\kV}_{\ps}$ & ${\objV}_{\ps}$ & $\eta$ & t & ${\kV}_{\fs}$ & ${\objV}_{\fs}$ & ${\kV}_{\ps}$ & ${\objV}_{\ps}$ & HIFOO & PENBMI\\
		\hline
		\hline
        AC1  & Inf     &1e0 & 0.19 & 76   & 0.038  & 131 & 0.034  & 1e0 & 0.18  & 76  & 0.038  & 131  & 0.034  & 1e0 & 0.15 & 102 & 0.037   & 152 & 0.034   &    -     &  0.360   \\ \hline
        AC2  & Inf     &1e0 & 0.17 & 76   & 0.038  & 131 & 0.034  & 1e0 & 0.19  & 76  & 0.038  & 131  & 0.034  & 1e0 & 0.17 & 102 & 0.037   & 152 & 0.034   & 0.050   &  0.300   \\ \hline
        AC4  & Inf     &1e4 & 0.21 & 1    & 11.026 &  1  & 11.026 & 1e4 & 0.20  &  1  & 11.026 &  1   & 11.026 & 1e4 & 0.17 & 1   & 11.026  & 1   & 11.026  &    -     &  11.014 \\ \hline
        AC6  & 24.606  &5e1 & 0.19 &  29  & 2.895  & 34  & 2.883  & 1e2 & 0.22  & 21  & 2.922  & 37   & 2.897  & 1e2 & 0.18 & 74  & 2.869   & 75  & 2.868   & 3.798   &  -      \\ \hline
        AC7  & Inf     &1e2 & 0.23 &  39  & 0.051  & 81  & 0.048  & 1e2 & 0.21  & 39  & 0.051  & 64   & 0.049  & 1e2 & 0.19 & 61  & 0.053   & 250 & 0.052   & 0.052   &  1.184       \\ \hline
        AC15 & 176.455 &1e0 & 0.18 &  27  & 2.554  & 37  & 1.908  & 1e0 & 0.18  & 33  & 2.580  & 43   & 1.896  & 2e0 & 0.17 & 25  & 2.669   & 76  & 1.776   & 12.612  &  353.728 \\ \hline
        AC17 & 10.265  &1e-1& 0.17 &  23  & 1.592  & 32  & 1.560  & 1e0 & 0.17  & 10  & 2.004  & 16   &   1.66 & 1e4 & 0.18 & 6   & 2.246   & 6   & 2.246   & 12.298  &  1.534  \\ \hline
        NN2  & Inf     &1e0 & 0.16 &   3  & 1.206  &  10 & 1.189  & 1e0 & 0.16  & 3   & 1.206  & 7    & 1.189  & 1e0 & 0.16 &  11 & 1.189   & 18  & 1.189   & 1.565   &  1.189  \\ \hline
        NN4  & 5.563   &1e1 & 0.17 &   8  & 2.062  &  33 & 1.928  & 1e1 & 0.17  & 9   & 2.062  & 34   & 1.926  & 5e1 & 0.16 &  10 & 2.159   & 47  & 1.964   & 1.875   &  1.832  \\ \hline
        NN8  & 5.922   & -  &   -  &  -   &    -   &  -  &    -   &  -  &   -   & -   &    -   &  -   &   -    & 1e1 & 0.16 & 5   & 1.772   & 36  & 1.596   & 2.279   &  1.510  \\ \hline
        NN11 & 0.142   &2e0 & 0.96 &   1  & 0.153  &  8  & 0.142  & 5e-1& 0.61  & 2   & 0.142  & 3    & 0.142  & 2e0 & 0.48 &  2  & 0.142   & 3   & 0.142   & 0.118   &  0.149  \\ \hline
        NN15 & Inf     &5e1 & 0.20 &  66  & 0.283  & 180 & 0.069  & 2e1 & 0.20  & 100 & 0.283  & 131  & 0.078  & 5e1 & 0.18 &  66 & 0.286   & 195 & 0.073   &  0.049  &  0.000  \\ \hline
        NN16 & Inf     &2e2 & 1.06 & 5    & 0.134  & 85  & 0.125  & 5e0 & 0.51  & 10  & 0.127  & 25   & 0.120  & 1e0 & 0.35 &  15 & 0.119   & 30  & 0.119   & 0.291   &  -       \\ \hline
        DIS1 & 5.149   &5e0 & 0.23 &  3   & 1.761  & 156 & 0.996  & 1e0 & 0.24  & 26  & 1.188  & 202  & 0.829  & 1e0 & 0.18 & 214 & 0.820   & 214 & 0.820   & 2.660   &  -       \\ \hline
        DIS3 & 11.653  &5e0 & 0.17 & 21   & 1.633  & 89  & 1.102  & 5e0 & 0.19  & 26  & 1.605  & 38   & 1.275  & 5e0 & 0.16 & 31  & 1.553   & 92  & 1.097   & 1.839   &  303.850  \\ \hline
        AGS  & 7.041   &1e3 & 0.38 & 217  & 7.057  & 226 & 7.057  & 1e3 & 0.47  & 217 &  7.056 & 227  & 7.055  &  -  &   -  & -   &     -   &  -  &    -    & 6.995   &  6.973  \\ \hline
        PSM  & 3.847   &2e-1& 0.18 & 23   & 0.242  & 250 & 0.072  & 2e-1& 0.18  & 24  & 0.237  & 250  & 0.072  & 5e-1& 0.15 &  10 & 0.496   & 250 & 0.091   & 1.503   &  0.004  \\ \hline
        BDT1 & 0.039   &5e-1& 0.27 & 74   & 0.016  & 250 & 0.005  & 5e-1& 0.24  & 78  & 0.015  & 250  & 0.005  & 1e0 & 0.20 & 206 & 0.013   & 250 & 0.009   & 0.010   &  0.006  \\ \hline
	\end{tabular}
    }
	\label{tab:ODC_Cent_Htwo}
\end{table*}

\begin{table*}[!t] 
\footnotesize
	\centering
	\caption{Results of centralized $\Hcal_{\infty}$ controller design for COMPl\textsubscript{e}ib models.}
     \scalebox{0.95}{
	\begin{tabular}{ |@{\;}c@{\;}|@{\;}c@{\;}
    |@{\;}|
    @{\!\!\:\,}c@{\;}|@{\;}c@{\;}|@{\;}c@{\;}|@{\;}c@{\;}|@{\;}c@{\;}|@{\;}c@{\;}
    |@{\;}|
    @{\;}c@{\;}|@{\;}c@{\;}|@{\;}c@{\;}|@{\;}c@{\;}|@{\;}c@{\;}|@{\;}c@{\;}
    |@{\;}|
    @{\;}c@{\;}|@{\;}c@{\;}|@{\;}c@{\;}|@{\;}c@{\;}|@{\;}c@{\;}|@{\;}c@{\;}
    |@{\;}|
    @{\;}c@{\;}|@{\;}c@{\;}|}
		\hline
		\multirow{2}{*}{Name} & \multirow{2}{*}{$\lVert\Gbfc_{\sm{\mathrm{ol}}{6}}{\rVert}_{\infty}$} & \multicolumn{6}{c|@{\;}|@{\;}}{SDP} & \multicolumn{6}{c|@{\;}|@{\;}}{SOCP} & \multicolumn{6}{c|@{\;}|@{\;}}{Parabolic} & \multicolumn{2}{c|}{Competitors}\\
		\cline{3-22}
 		& &  $\eta$ & t & ${\kV}_{\fs}$ & ${\objV}_{\fs}$ & ${\kV}_{\ps}$ & ${\objV}_{\ps}$ & $\eta$ & t & ${\kV}_{\fs}$ & ${\objV}_{\fs}$ & ${\kV}_{\ps}$ & ${\objV}_{\ps}$ & $\eta$ & t & ${\kV}_{\fs}$ & ${\objV}_{\fs}$ & ${\kV}_{\ps}$ & ${\objV}_{\ps}$ & HIFOO & PENBMI\\
		\hline
		\hline
        AC1  & 2.167   &5e-1& 0.19 & 14 & 0.000 & 250 & 0.000  & 5e-1& 0.21  & 15 & 0.000 & 126 & 0.000  & 5e-1 & 0.17 & 14 & 0.000 & 110 & 0.000   & 0.000  &   0.008 \\ \hline
        AC2  & 2.167   &1e3 & 0.18 & 83 & 0.602 & 250 & 0.296  & 5e2 & 0.21  & 94 & 0.440 & 250 & 0.238  & 5e2  & 0.16 & 250& 0.237 & 250 & 0.237   & 0.111  &  0.118  \\ \hline
        AC4  & 69.990  &1e0 & 0.18 & 2  & 70.078& 8   & 69.990 & 1e0 & 0.19  & 2  & 70.078& 4   & 69.990 & 1e0  & 0.18 & 2  & 70.078& 11  & 69.990  & 0.935  &   -      \\ \hline
        AC6  & 391.782 &5e2 & 0.21 & 174& 6.584 & 250 & 4.410  & 5e2 & 0.21  & 169& 6.636 & 250 & 4.385  & -    &   -  &  - &    -   &  -  &    -     & 4.113  &  4.113  \\ \hline
        AC7  & 0.042   &1e0 & 0.22 & 20 & 0.000 & 110 & 0.000  & 2e0 & 0.20  & 32 & 0.000 & 41  & 0.000  & 1e0  & 0.19 & 21 & 0.000 & 173 & 0.000   & 0.064  &   -      \\ \hline
        AC15 & 2.4e3   &5e2 & 0.20 & 165& 39.945& 250 & 20.250 & 5e2 & 0.23  & 166& 39.983& 250 & 20.283 & -    &   -  & -  &    -   &  -  &    -     & 16.865 &  15.168 \\ \hline
        AC17 & 30.823  &1e4 & 0.19 & 26 & 79.491& 250 & 7.748  & 1e4 & 0.19  & 26 & 79.682& 250 & 7.748  & 1e4  & 0.19 &  34& 71.867& 250 &  7.640  & 16.639 &  14.855 \\ \hline
        NN2  & Inf     &1e0 & 0.17 & 17 & 2.224 & 47  & 2.221  & 1e0 & 0.15  & 17 & 2.224 & 56  & 2.221  & 5e0  & 0.15 & 5  & 2.769 & 16  & 2.222   & 2.220  &  2.221  \\ \hline
        NN4  & 31.043  &5e0 & 0.17 & 18 & 1.792 & 67  & 1.416  & 5e0 & 0.17  & 20 & 1.729 & 72  & 1.412  & 5e0  & 0.16 & 32 & 1.551 & 75  & 1.411   & 1.369  &  1.358  \\ \hline
        NN8  & 46.508  &5e3 & 0.18 & 6  & 59.946& 250 & 3.585  & 5e3 & 0.17  & 6  & 59.910& 250 & 3.585  & 5e3  & 0.17 & 12 & 43.637& 250 & 3.587   & 3.387  &   -      \\ \hline
        NN11 & 0.170   &5e0 & 2.26 & 48 & 0.460 & 180 & 0.169  & 5e0 & 1.15  & 50 & 0.440 & 250 & 0.190  & 5e0  & 1.37 & 179& 0.193 & 207 & 0.169   & 0.107  &  0.124  \\ \hline
        NN15 & Inf     &1e-2& 0.19 & 109& 0.278 & 250 & 0.130  & 1e-2& 0.23  & 108& 0.284 & 250 & 0.130  & 1e-2 & 0.17 & 152& 0.232 & 250 & 0.134   & 0.098  &  0.098  \\ \hline
        NN16 & 6.4e14  &1e3 & 0.88 & 56 & 0.575 & 81  & 0.559  & 1e3 & 0.53  & 77 & 0.566 & 95  & 0.557  &  -   &   -  &  - &    -   &  -  &    -     & 1.012  &   -      \\ \hline
        DIS1 & 17.320  &1e1 & 0.45 & 54 & 4.574 & 91  & 4.286  & 1e1 & 0.41  & 52 & 4.580 & 90  & 4.287  & 2e1  & 0.29 & 99 & 4.626 & 104 & 4.564   & 4.182  &   -      \\ \hline
        DIS3 & 32.069  &1e1 & 0.30 & 73 & 2.613 & 116 & 1.302  & 1e1 & 0.25  & 74 & 2.625 & 101 & 1.308  & 2e2  & 0.19 & 37 & 5.876 & 215 & 1.350   & 1.341  &  1.275  \\ \hline
        AGS  & 8.182   &2e3 & 0.35 & 144& 8.872 & 188 & 8.192  & 2e3 & 0.34  & 150& 8.751 & 194 & 8.189  & 1e4  & 0.27 & 209& 10.765& 250 & 9.486   & 8.173  &  8.173  \\ \hline
        PSM  & 4.232   &5e0 & 0.17 & 20 & 1.185 & 41  & 0.921  & 5e0 & 0.17  & 20 & 1.195 & 41  & 0.921  & 1e1  & 0.16 & 25 & 1.231 & 52  & 0.921   & 0.920  &  0.920  \\ \hline
        BDT1 & 5.142   &1e-2& 0.84 & 77 & 0.565 & 250 & 0.311  & 1e0 & 0.87  & 22 & 5.649 & 250 & 0.873  & 5e-2 & 0.53 & 202& 0.523 & 250 & 0.431   & 0.266  &  0.266  \\ \hline
	\end{tabular}
   }
	\label{tab:ODC_Cent_Hinf}
\end{table*}

\subsection{Case Study I: Centralized Controller}
We use the proposed sequential method to find unstructured static output-feedback controllers that stabilize the plant and minimize the norm of the closed-loop system. This controller is allowed to use the entire measurements to generate the control decisions. Numerical results for $\Hcal_2$ and $\Hcal_{\infty}$ controller design problems are reported in \cref{tab:ODC_Cent_Htwo,tab:ODC_Cent_Hinf}, respectively. 

The first two columns of the tables contain the model names and the corresponding norm of the open-loop systems. We compute the open-loop norms based on the following system:
\begin{equation}\nonumber
\begin{aligned}
\Gbfc_{\mathrm{ol}}
=\begin{bmatrix} \Abfc & \Bbfc_{\sm{1}{6}}\\ \Cbfc_{\sm{1}{6}} & \Dbfc_{\sm{11}{6}}
\end{bmatrix}.
\end{aligned}
\end{equation}
The subsequent sections in Tables \ref{tab:ODC_Cent_Htwo} and \ref{tab:ODC_Cent_Hinf} consist of four different sections labeled as SDP, SOCP, parabolic, and competitors. The first three sections show the numerical results of the proposed method equipped with different relaxations. The last section contains the numerical results of HIFOO and PENBMI. The following numbers are reported in the SDP, SOCP and parabolic sections:
\begin{itemize}
\item $\eta$ denotes the choice of penalty parameter in \cref{eq:ODC_GenRelaxation_Pen_obj}. This parameter is chosen from the set $\{1\times 10^i, 2\times 10^i, 5\times 10^i\}_{\sm{i=-2}{5.5}}^{\sm{4}{5.5}}$ in all experiments. 
\item $t$ denotes the average run time of solving each round of penalized convex relaxation in \cref{al:ODC_alg_01}.
\item ${\kV}_{\fs}$ denotes the number of rounds necessary to obtain a feasible solution for the original BMI (i.e., the first round whose resulting solution satisfies $\mathbf{X}=\mathbf{x}\mathbf{x}^{\top}$) and ${\objV}_{\fs}$ represents the corresponding objective value at round ${\kV}_{\fs}$ (without the penalty term).
\item ${\kV}_{\ps}$ and ${\objV}_{\ps}$, respectively, denote the round number at which the stopping criteria is met and the corresponding objective value.
\end{itemize}
In all of the experiments, we terminate the sequential penalized relaxation when the percentage objective value improvement between two consecutive rounds is less than 0.1 for $\Hcal_2$ and 0.05 for $\Hcal_{\infty}$, or if the number of rounds exceeds 250. For cases where $\Bbfc_{\sm{1}{6}}\Bbfc_{\sm{1}{6}}^{\!\top}$ is not positive definite, we substitute it with the matrix $\Bbfc_{\sm{1}{6}}\Bbfc_{\sm{1}{6}}^{\!\top} + 10^{-5}\times\Ibf$. 

Given an optimal solution $\accentset{\ast}{\xbf}$ for either \cref{eq:ODC_Htwo_obj,eq:ODC_Htwo_con_01} or \cref{eq:ODC_Hinf_obj,eq:ODC_Hinf_con_01}, let $\accentset{\ast}{\hbfc}$ represent the entries of $\accentset{\ast}{\xbf}$ corresponding to the controller element. We consider $\Kbfc(\msh\accentset{\ast}{\hbfc}\msh)$ as a stabilizing controller for the plant $\Gbfc$ if the real part of all eigenvalues of $\Abfc+\Bbfc\Kbfc(\msh\accentset{\ast}{\hbfc}\msh)\Cbfc$ are smaller than $10^{-5}$.
\subsection{Case Study II: Decentralized Controller}
This case study is concerned with the design of decentralized controllers. Despite the centralized case, this controller only have access to a subset of measurements to generate the control commands. In this experiment, we only consider those models in which the control commands vector $\ubfc$ and the sensor measurements vector $\ybfc$ are of the same dimensions. We apply \cref{al:ODC_alg_01} to find $\Hcal_2$ and $\Hcal_{\infty}$ static output-feedback controllers with diagonal patterns. The results of this experiment are reported in \cref{tab:ODC_DeCent_Htwo,tab:ODC_DeCent_Hinf}. 

As the results indicate, the performance of the proposed sequential scheme equipped with SDP, SOCP, and parabolic relaxations provide promising performance for both centralized and decentralized cases compared to both PENBMI and HIFOO packages (smaller norm means better performance). 

\begin{table*}[!t] 
\footnotesize
	\centering
	\caption{Results of fully decentralized $\Hcal_2$ controller design for COMPl\textsubscript{e}ib models.}
	\begin{tabular}{ |@{\;}c@{\;}|@{\;}c@{\;}
    |@{\;}|
    @{\!\!\:\,}c@{\;}|@{\;}c@{\;}|@{\;}c@{\;}|@{\;}c@{\;}|@{\;}c@{\;}|@{\;}c@{\;}
    |@{\;}|
    @{\;}c@{\;}|@{\;}c@{\;}|@{\;}c@{\;}|@{\;}c@{\;}|@{\;}c@{\;}|@{\;}c@{\;}
    |@{\;}|
    @{\;}c@{\;}|@{\;}c@{\;}|@{\;}c@{\;}|@{\;}c@{\;}|@{\;}c@{\;}|@{\;}c@{\;}
    |@{\;}|
    @{\;}c@{\;}|@{\;}c@{\;}|}
		\hline
		\multirow{2}{*}{Name} & \multirow{2}{*}{$\lVert\Gbfc_{\sm{\mathrm{ol}}{6}}{\rVert}_{2}$} & \multicolumn{6}{c|@{\;}|@{\;}}{SDP} & \multicolumn{6}{c|@{\;}|@{\;}}{SOCP} & \multicolumn{6}{c|@{\;}|@{\;}}{Parabolic} & \multicolumn{2}{c|}{Competitors}\\
		\cline{3-22}
 		& &  $\eta$ & t & ${\kV}_{\fs}$ & ${\objV}_{\fs}$ & ${\kV}_{\ps}$ & ${\objV}_{\ps}$ & $\eta$ & t & ${\kV}_{\fs}$ & ${\objV}_{\fs}$ & ${\kV}_{\ps}$ & ${\objV}_{\ps}$ & $\eta$ & t & ${\kV}_{\fs}$ & ${\objV}_{\fs}$ & ${\kV}_{\ps}$ & ${\objV}_{\ps}$ & HIFOO & PENBMI\\
		\hline
		\hline
        AC1  & Inf    & -   &   -  & -  &   -   & -   & -       & -   &   -   & -  &    -  & -   &   -    & 1e1 & 0.14  & 59 & 0.084 & 196 & 0.047   & 0.054  &  0.039  \\ \hline
        AC2  & Inf    & -   &   -  & -  &   -   & -   & -       & -   &   -   & -  &    -  & -   &   -    & 1e1 & 0.15  & 59 & 0.084 & 196 & 0.047   & 0.090  &  0.039  \\ \hline
        NN2  & Inf    & 2e0 & 0.13 & 2  & 1.220 & 8   & 1.189   & 2e0 & 0.14  & 2  & 1.220 & 8   & 1.189  & 1e0 & 0.13  & 11 & 1.191 & 18  & 1.189   & 1.565  &  1.189  \\ \hline
        NN8  & 5.9220 & 5e0 & 0.14 & 3  & 1.864 & 10  & 1.839   & 5e0 & 0.14  & 3  & 1.864 & 10  & 1.839  & 5e0 & 0.14  & 4  & 1.864 & 10  & 1.840   & 2.365  &  1.838  \\ \hline
        NN15 & Inf    & 1e0 & 0.16 & 238& 0.122 & 250 & 0.082   & 2e0 & 0.16  & 229& 0.275 & 250 & 0.077  & 5e1 & 0.15  & 67 & 0.283 & 183 & 0.070   & 0.049  &  0.000  \\ \hline
        NN16 & Inf    & 1e0 & 0.19 & 5  & 0.139 & 54  & 0.121   & 2e1 & 0.21  & 6  & 0.140 & 66  & 0.122  & 5e-1& 0.16  & 13 & 0.132 & 75  & 0.119   & 0.488  &  0.119  \\ \hline
        DIS1 & 5.1491 & 1e2 & 0.19 & 9  & 2.111 & 54  & 1.783   & 1e2 & 0.20  & 9  & 2.115 & 56  & 1.781  & 1e3 & 0.18  & 21 & 2.351 & 132 & 1.874   & 2.991  &6982.151 \\ \hline
        DIS2 & Inf    & 5e0 & 0.15 & 24 & 1.501 & 119 & 0.512   & 5e0 & 0.15  & 24 & 1.501 & 119 & 0.512  & 5e0 & 0.14  & 47 & 1.368 & 120 & 0.487   & 2.047  &  0.377  \\ \hline
        DIS3 & 11.6538& 1e1 & 0.16 & 20 & 2.039 & 107 & 1.385   & 1e1 & 0.17  & 20 & 2.038 & 107 & 1.385  & 1e1 & 0.15  & 43 & 1.976 & 95  & 1.370   & 2.286  &   -      \\ \hline
        AGS  & 7.0412 & 2e3 & 0.29 & 187& 7.149 & 187 & 7.149   & 1e3 & 0.34  & 214& 7.038 & 214 & 7.038  &  -  &   -    & -  &   -    &  -  &     -    & 7.029  &  7.032  \\ \hline
        BDT1 & 0.0397 & 1e0 & 0.24 &  34& 0.026 & 250 & 0.007   & 5e-1& 0.21  & 85 & 0.018 & 250 & 0.006  & 2e0 & 0.19  &119 & 0.020 & 250 & 0.008   & 0.010  & 0.000      \\ \hline
	\end{tabular}
	\label{tab:ODC_DeCent_Htwo}
\end{table*}

\begin{table*}[!t] 
\vspace{0.11cm}
\footnotesize
	\centering
	\caption{Results of fully decentralized $\Hcal_{\infty}$ controller design for COMPl\textsubscript{e}ib models.}
	\begin{tabular}{ |@{\;}c@{\;}|@{\;}c@{\;}|@{\;}|@{\!\!\:\,}c@{\;}|@{\;}c@{\;}|@{\;}c@{\;}|@{\;}c@{\;}|@{\;}c@{\;}|@{\;}c@{\;}|@{\;}|@{\;}c@{\;}|@{\;}c@{\;}|@{\;}c@{\;}|@{\;}c@{\;}|@{\;}c@{\;}|@{\;}c@{\;}|@{\;}|@{\;}c@{\;}|@{\;}c@{\;}|@{\;}c@{\;}|@{\;}c@{\;}|@{\;}c@{\;}|@{\;}c@{\;}|@{\;}|@{\;}c@{\;}|@{\;}c@{\;}|}
		\hline
		\multirow{2}{*}{Name} & \multirow{2}{*}{$\lVert\Gbfc_{\sm{\mathrm{ol}}{6}}{\rVert}_{\infty}$} & \multicolumn{6}{c|@{\;}|@{\;}}{SDP} & \multicolumn{6}{c|@{\;}|@{\;}}{SOCP} & \multicolumn{6}{c|@{\;}|@{\;}}{Parabolic} & \multicolumn{2}{c|}{Competitors}\\
		\cline{3-22}
 		& &  $\eta$ & t & ${\kV}_{\fs}$ & ${\objV}_{\fs}$ & ${\kV}_{\ps}$ & ${\objV}_{\ps}$ & $\eta$ & t & ${\kV}_{\fs}$ & ${\objV}_{\fs}$ & ${\kV}_{\ps}$ & ${\objV}_{\ps}$ & $\eta$ & t & ${\kV}_{\fs}$ & ${\objV}_{\fs}$ & ${\kV}_{\ps}$ & ${\objV}_{\ps}$ & HIFOO & PENBMI\\
		\hline
		\hline
        AC1  & 2.167  & 5e0 & 0.14 & 37 & 0.164 & 250 & 0.056  & 5e0 & 0.15  & 37 & 0.167 & 250 & 0.057  & 1e1  & 0.15 & 55 & 0.162 & 250 & 0.064    & 0.067  &  0.014   \\ \hline
        AC2  & 2.167  & 5e2 & 0.15 & 160& 0.422 & 250 & 0.319  & 5e2 & 0.16  & 143& 0.453 & 250 & 0.319  & 5e3  & 0.14 & 149& 0.921 & 250 & 0.661    & 0.661  &  0.167   \\ \hline
        NN2  & Inf    & 5e0 & 0.13 & 5  & 2.758 & 16  & 2.222  & 5e0 & 0.13  & 5  & 2.758 & 16  & 2.222  & 5e0  & 0.13 & 5  & 2.769 & 16  & 2.222    &2.220   &  2.221   \\ \hline
        NN8  & 46.508 & 1e3 & 0.14 & 6  & 30.859& 194 & 3.405  & 1e3 & 0.14  & 6  & 30.859& 194 & 3.405  & 2e3  & 0.14 & 11 & 30.752& 236 & 3.461    &3.272   &  3.746   \\ \hline
        NN15 & Inf    & 1e-2& 0.14 & 107& 0.288 & 250 & 0.130  & 1e-2& 0.15  & 107& 0.288 & 250 & 0.130  & 1e-2 & 0.14 & 152& 0.231 & 250 & 0.135    &0.100   &  0.100   \\ \hline
        NN16 & 6.4e14 & 2e3 & 0.20 & 35 & 1.143 & 59  & 0.959  & 5e3 & 0.20  & 54 & 1.342 & 89  & 0.959  & 5e3  & 0.16 & 78 & 1.118 & 111 & 0.960    &0.956   &  0.957   \\ \hline
        DIS1 & 17.320 & 1e2 & 0.18 & 30 & 11.591& 64  & 7.175  & 5e1 & 0.21  & 43 & 10.235& 67  & 7.186  & 1e2  & 0.14 & 52 & 10.028& 84  & 7.223    & 7.165  & 6.843    \\ \hline
        DIS3 & 32.069 & 5e1 & 0.16 & 63 & 5.831 & 115 & 1.689  & 5e1 & 0.17  & 80 & 4.996 & 115 & 1.696  & 5e1  & 0.14 & 37 & 4.510 & 91  & 1.690    & 1.655  & 1.656    \\ \hline
        AGS  & 8.182  & 5e3 & 0.26 & 106& 13.143& 231 & 8.209  & 2e3 & 0.25  & 148& 8.697 & 191 & 8.196  & 1e3  & 0.20 & 202& 11.041& 250 & 9.420    & 8.173  & 8.173    \\ \hline
        BDT1 & 5.142  & 2e-1& 0.28 & 30 & 2.704 & 250 & 0.575  & 2e-1& 0.31  & 30 & 2.704 & 250 & 0.575  & 1e-1 & 0.23 & 168& 0.737 & 250 & 0.551    & 0.266  & 0.266    \\ \hline
	\end{tabular}
	\label{tab:ODC_DeCent_Hinf}
\end{table*}


\subsection{Case Study III: Choice of Penalty Parameter $\eta$}
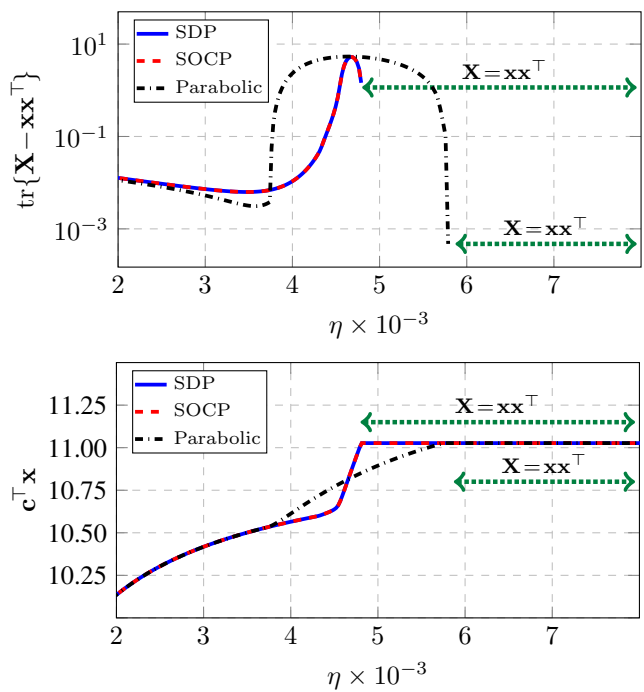
\begin{figure}
\captionsetup[subfigure]{position=b}
\centering
\begin{tikzpicture}[scale=1]
\centering
\begin{axis}[
		width = 0.48\textwidth,
		height = 0.28\textwidth,
 		xmode=normal,
		ymode=log,	
    	ylabel={$\tr\{\Xbf\!-\!\xbf\xbf^{\top}\!\}$},
    	ylabel shift =0cm,
 		xmin=2000,
	 	xmax=8000,
        xtickmax = 7500,
	 	ymax= 50,
        xtick scale label code/.code={\pgfmathparse{int(#1)}$\eta \times 10^{-\pgfmathresult}$},
 		every x tick scale label/.style={at={(xticklabel cs:0.5)}, anchor = north},
		scaled x ticks = base 10:-3,
    	legend pos=north west,
    	ymajorgrids=true,
		xmajorgrids=true,
    	grid style=dashed,
		legend cell align={left},
        legend image post style={scale=0.7},
        legend style={font=\scriptsize, inner xsep=0.2pt, inner ysep=0.1pt},
]

\draw[<->, densely dotted, line width= 0.06cm, Fcolor] (axis cs:4800.22,1.15)--(axis cs:7950,1.15) node [black, pos=0.52, above] {\footnotesize \raisebox{-0.25ex}[0pt][-0.25ex]{$\Xbf\!=\!\xbf\xbf^{\top}$}};
\draw[<->, densely dotted, line width= 0.06cm, Fcolor] (axis cs:5870.92,0.000473)--(axis cs:7950,0.000473) node [black, pos=0.5, above] {\footnotesize \raisebox{-0.25ex}[0pt][-0.25ex]{$\Xbf\!=\!\xbf\xbf^{\top}$}};
\addplot[color=blue,line width=0.5mm] table[x=X,y=FSOCP,col sep=comma] {AC4.csv};
\addplot[dashed,color=red,line width=0.5mm] table[x=X,y=FSDP,col sep=comma] {AC4.csv};
\addplot[dash dot,color=black,line width=0.5mm] table[x=X,y=FPARA,col sep=comma] {AC4.csv};

\legend{$\mathrm{SDP}$,$\mathrm{SOCP}$,$\mathrm{Parabolic}$}
\end{axis}
\label{plt:plot_feasibility_violation}
\end{tikzpicture}

\vspace{0.2cm}

\begin{tikzpicture}[scale=1]
\centering
\begin{axis}[
		width = 0.48\textwidth,
		height = 0.28\textwidth,
 		xmode=normal,
		ymode=normal,	
        ylabel={$\cbf^{\!\top}\xbf$},
		ylabel shift = 0cm,
 		xmin=2000,
  		xmax=8000,
    	ymin = 10.0,
     	ymax = 11.5,
        ytickmax = 11.4,,
        xtickmax = 7500,
        xtick scale label code/.code={\pgfmathparse{int(#1)}$\eta \times 10^{-\pgfmathresult}$},
 		every x tick scale label/.style={at={(xticklabel cs:0.5)}, anchor = north},
		scaled x ticks = base 10:-3,
		ytick={10.00, 10.25, 10.50, 10.75, 11.00, 11.25, 11.50},
    	yticklabels={~, 10.25, 10.50, 10.75, 11.00, 11.25, 11.50},
    	legend pos= north west,
    	ymajorgrids=true,
		xmajorgrids=true,
    	grid style=dashed,
		legend cell align={left},
        legend image post style={scale=0.7},
        legend style={font=\scriptsize, inner xsep=0.2pt, inner ysep=0.1pt},
]
\draw[<->, densely dotted, line width= 0.06cm, Fcolor] (axis cs:4800.22,11.15)--(axis cs:7950,11.15) node [black, pos=0.5, above] {\footnotesize \raisebox{-0.2ex}[0pt][-0.2ex]{$\Xbf\!=\!\xbf\xbf^{\top}$}};
\draw[<->, densely dotted, line width= 0.06cm, Fcolor] (axis cs:5870.92,10.8)--(axis cs:7950,10.8) node [black, pos=0.5, above] {\footnotesize \raisebox{-0.2ex}[0pt][-0.2ex]{$\Xbf\!=\!\xbf\xbf^{\top}$}};
\addplot[color=blue,line width=0.5mm] table[x=X,y=CSOCP,col sep=comma] {AC4.csv};
\addplot[dashed,color=red,line width=0.5mm] table[x=X,y=CSDP,col sep=comma] {AC4.csv};
\addplot[dash dot,color=black,line width=0.5mm] table[x=X,y=CPARA,col sep=comma] {AC4.csv};

\legend{$\mathrm{SDP}$,$\mathrm{SOCP}$,$\mathrm{Parabolic}$}
\end{axis}
\end{tikzpicture}
\caption{Sensitivity of the penalized convex relaxations to the choice of penalty parameter $\eta$ ($\Hcal_2$ centralized controller design for ``AC4'' model).}
\label{plt:ODC_AC4}
\end{figure}

\begin{figure}
\captionsetup[subfigure]{position=b}
\centering
\begin{tikzpicture}[scale=1]
\centering
\begin{axis}[
		width = 0.48\textwidth,
		height = 0.28\textwidth,
 		xmode=normal,
		ymode=log,	
    	ylabel={$\tr\{\Xbf-\xbf\xbf^{\top}\}$},
    	ylabel shift = 0cm,
  		xmin=0,
 	 	xmax=100,
        xtick scale label code/.code={\pgfmathparse{int(#1)}$\eta \times 10^{-\pgfmathresult}$},
 		every x tick scale label/.style={at={(xticklabel cs:0.5)}, anchor = north},
		xtickmax = 90,
		scaled x ticks = base 10:-1,
		legend pos = north east,
    	ymajorgrids=true,
		xmajorgrids=true,
    	grid style=dashed,
		legend cell align={left},
        legend image post style={scale=0.7},
        legend style={font=\scriptsize, inner xsep=0.2pt, inner ysep=0.1pt},
]
\draw[<->, densely dotted, line width= 0.06cm, Fcolor] (axis cs:7.5,0.00001)--(axis cs:99,0.00001) node [black, pos=0.58, above] {\footnotesize\raisebox{-0.2ex}[0pt][-0.2ex]{$\Xbf\!=\!\xbf\xbf^{\top}$}};
\draw[<->, densely dotted, line width= 0.06cm, Fcolor] (axis cs:48,0.000001)--(axis cs:99,0.000001) node [black, pos=0.5, above] {\footnotesize \raisebox{-0.2ex}[0pt][-0.2ex]{$\Xbf\!=\!\xbf\xbf^{\top}$}};
\addplot[color=blue,line width=0.5mm] table[x=X,y=FSOCP,col sep=comma] {BDT1.csv};
\addplot[dashed,color=red,line width=0.5mm] table[x=X,y=FSDP,col sep=comma] {BDT1.csv};
\addplot[dash dot,color=black,line width=0.5mm] table[x=X,y=FPARA,col sep=comma] {BDT1.csv};

\legend{$\mathrm{SDP}$,$\mathrm{SOCP}$,$\mathrm{Parabolic}$}
\end{axis}
\label{plt:plot_feasibility_violation}
\end{tikzpicture}

\vspace{0.2cm}

\begin{tikzpicture}[scale=1]
\centering
\begin{axis}[
		width = 0.48\textwidth,
		height = 0.28\textwidth,
 		xmode=normal,
		ymode=normal,	
        ylabel={$\cbf^{\!\top}\xbf$},
		ylabel shift = 0cm,
  		xmin = 0,
   		xmax = 100,
     	ymin =0,
      	ymax = 0.085,
        xtick scale label code/.code={\pgfmathparse{int(#1)}$\eta \times 10^{-\pgfmathresult}$},
 		every x tick scale label/.style={at={(xticklabel cs:0.5)}, anchor = north},
		xtickmax = 90,
		scaled x ticks = base 10:-1,
        scaled y ticks = false,
        tick label style={/pgf/number format/fixed},
    	legend pos = north east,
    	ymajorgrids = true,
		xmajorgrids = true,
    	grid style = dashed,
		legend cell align = {left},
        legend image post style = {scale=0.7},
        legend style={font=\scriptsize, inner xsep=0.2pt, inner ysep=0.1pt},
]
\draw[<->, densely dotted, line width= 0.06cm, Fcolor] (axis cs:7.4,0.047)--(axis cs:99,0.047) node [black, pos=0.5, above] {\footnotesize \raisebox{-0.25ex}[0pt][-0.25ex]{$\Xbf\!=\!\xbf\xbf^{\top}$}};
\draw[<->, densely dotted, line width= 0.06cm, Fcolor] (axis cs:48,0.0257)--(axis cs:99,0.0257) node [black, pos=0.5, above] {\footnotesize \raisebox{-0.25ex}[0pt][-0.25ex]{$\Xbf\!=\!\xbf\xbf^{\top}$}};

\addplot[color=blue,line width=0.5mm] table[x=X,y=CSOCP,col sep=comma] {BDT1.csv};
\addplot[dashed,color=red,line width=0.5mm] table[x=X,y=CSDP,col sep=comma] {BDT1.csv};
\addplot[dash dot,color=black,line width=0.5mm] table[x=X,y=CPARA,col sep=comma] {BDT1.csv};

\legend{$\mathrm{SDP}$,$\mathrm{SOCP}$,$\mathrm{Parabolic}$}
\end{axis}
\end{tikzpicture}
\caption{Sensitivity of the penalized convex relaxations to the choice of penalty parameter $\eta$ ($\Hcal_2$ centralized controller design for ``BDT1'' model).}
\label{plt:ODC_BDT1}
\end{figure}
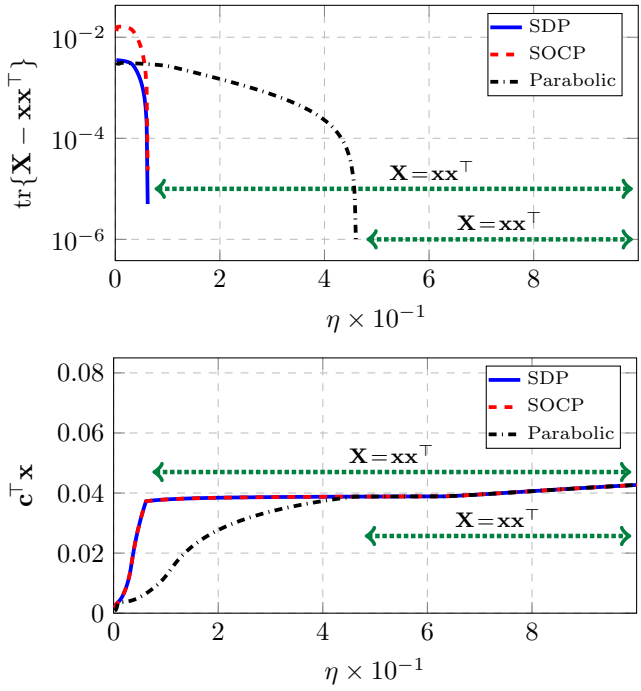

This case study investigates the sensitivity of different convex relaxations to the choice of regularization parameter $\eta$. To this end, one round of the penalized relaxation problem \cref{eq:ODC_GenRelaxation_Pen_obj,eq:ODC_GenRelaxation_Pen_con_02} (with zero initialization) is solved for a wide range of $\eta$ values. This experiment is run on the $\Hcal_2$-norm static output-feedback controller problem for the two models ''AC4'' and ``BDT1''. \cref{plt:ODC_AC4,plt:ODC_BDT1} show that if $\eta$ is small, none of the proposed penalized relaxations methods provide feasible solutions for BMI. As the value of $\eta$ increases, the feasibility violation $\mathrm{tr}\{\mathbf{X}-\mathbf{x}\mathbf{x}^{\top}\}$ abruptly vanishes once crossing a certain threshold. In our experiments, this limit has been close for SDP and SOCP relaxations, which is smaller than that of the parabolic relaxation. According to \cref{plt:ODC_AC4,plt:ODC_BDT1}, all three methods produce feasible points for a wide range of $\eta$ values, and within that range, the objective cost is not very sensitive to the choice of $\eta$. 
\section{Conclusion}
In this paper, we proposed a feasibility preserving sequential algorithm to solve optimization problems with linear objectives subject to bilinear matrix inequality (BMI) constraints. The proposed method can start from any arbitrary initial point to obtain feasible and near-optimal solutions. The performance of the proposed sequential method is tested on the problems of $\Hcal_2$ and $\Hcal_{\infty}$ control design for benchmark plants from COMPl\textsubscript{e}ib \cite{leibfritz2006compleib}. The numerical results verify the promising performance of our sequential penalized relaxation in comparison with the HIFOO and the PENBMI packages.
\appendices

\bibliographystyle{IEEEtran}
\bibliography{IEEEabrv,egbib}
\end{document}